\documentclass[leqno,11pt]{article}
\usepackage[spanish]{babel}
\usepackage[ansinew]{inputenc}
\usepackage{amscd}
\usepackage{amsthm}
\usepackage{amssymb}
\usepackage{amsmath}
\usepackage{graphics}
\usepackage{verbatim}
\usepackage{color}
\newtheorem{mydef}{Definition}
\newtheorem{myexa}{Example}
\newtheorem{mycexa}{Counterexample}
\newtheorem{mytheo}{Theorem}

\newtheorem*{myrem}{{\em Remark}}

\newcommand{\cind}{\perp\hspace{-1.25ex}\perp}

\newcommand{\pt}{\mbox{$\succ$\hspace{-1ex}$\longrightarrow$}}

\decimalpoint

\begin{document}

\begin{center}
	{\sc A Note on the Relationship Between Conditional and Unconditional Independence, and its Extensions for Markov Kernels}\vspace{2ex}\\
	A.G. Nogales and P. Pérez\vspace{2ex}\\
	Dpto. de Matem\'aticas, Universidad de Extremadura\\
	Avda. de Elvas, s/n, 06006--Badajoz, SPAIN.\\
	e-mail: nogales@unex.es
\end{center}
\vspace{.4cm}
\begin{quote}
	\hspace{\parindent} {\small {\sc Abstract.} Two known results on the relationship between conditional and unconditional independence are obtained as a consequence of the main result of this paper, a theorem that uses independence of Markov kernels  to obtain a minimal condition which added to conditional independence implies independence.		
		Some {  examples,} counterexamples and representation results are provided to clarify the concepts introduced and the propositions of the statement of the main theorem. Moreover, conditional independence and the mentioned results are extended to the framework of Markov kernels. }
\end{quote}

\vfill
\begin{itemize}
	\item[] \hspace*{-1cm} {\em AMS Subject Class.} (2010): {\em Primary\/} 60Exx
	{\em Secondary\/} 60J35
	\item[] \hspace*{-1cm} {\em Key words and phrases:} conditional independence, Markov kernel.

\end{itemize}

\newpage

\section{Introduction and basic definitions}

Conditional independence is a classical and familiar basic tool of both probability theory (think on Markov chains theory, for example) and mathematical statistics.  {  See, for instance, Dawid (1979) and Florens et al. (1990), where an extensive use of conditional independence is made in order to unify many seemingly unrelated concepts of statistical inference, either from the Bayesian and the frequentist point of views.  	
The introduction sections of Dawid (1979), Phillips (1988) and van Putten et al. (1985) list some of the main fields of application of the conditional independence relation: e.g.  econometric distribution theory, asymptotic studies of regression with non-ergodic processes,  the definition of a stochastic dynamic system and the stochastic realization problem, or, in a statistical framework, the areas of sufficiency, ancillarity, identification or invariance, among others. van Putten et al. (1985) includes also, among other interesting results, a systematic study of invariance properties of conditional independence under enlargement or reduction of the involved $\sigma$-fields, which keep some connection with the main problem raised in this paper. } 

It is well known that conditional independence does not imply, and it is not implied by, independence. We shall write $X\cind Y$ and $X\cind Y|Z$ for the independence of the random variables $X$ and $Y$ and its conditional independence given a third random variable $Z$, respectively. 

Section \ref{ci} contains the main result of this paper, Theorem \ref{teo2}, that uses independence of Markov kernels, a concept introduced by Nogales (2013a),  to obtain a minimal condition which added to conditional independence implies independence. 

This way the result becomes an improvement of two known results on the relationship between conditional and unconditional independence: one that constitutes the main goal of Phillips (1988), and another that is obtained as an immediate consequence of Theorem 2.2.10 of Florens et al. (1990) (or the Lemma 4.3 of Dawid (1979)), as it is remarked in Section \ref{cexa}. In this section some examples and counterexamples are also given to delimit the relations between the  three propositions of Theorem \ref{teo2}. 

 In this paper (Section \ref{extmk}) we also attack the problem of constructing a rigorous general theory of conditional independence in terms of Markov kernels; notice that Markov kernels are extensions of the concepts of both random variable and $\sigma$-field, and Theorem \ref{teo2} is here extended to this new framework. Dawid (1980) constructs a theory of conditional independence for ``statistical operations'', which is presented as a slight generalization of Markovian operator, which appears itself as a generalization of Markov kernel. Although this article also runs in the field of specialized mathematics, we hope the reader can find the development of conditional independence in the less abstract frame of Markov kernels (or transition probabilities) useful.
 
 A more general result than Theorem \ref{teo2} in terms of random variables is finally presented in Section \ref{last}. The introduced definition of conditional independence between Markov kernels is used to obtain a minimal condition which added to conditional independence of $X_1$ and $X_2$ given $X_3$ implies the conditional independence of $X_1$ and $X_2$ given $X_4$, provided $X_4$ is a function of $X_3$.

The paper is completed with some understandable reformulations of  several of the propositions considered. With the same purpose, some representation results of the introduced definitions for Markov kernels in terms of random variables are also facilitated. 

For ease of reading, the demonstrations will appear in a final section

In what follows  $(\Omega,\mathcal A)$, $(\Omega_1,\mathcal A_1)$,
and so on, will denote measurable spaces. A random variable is a map
$X:(\Omega,\mathcal A)\rightarrow (\Omega_1,\mathcal A_1)$ such that
$X^{-1}(A_1)\in\mathcal A$, for all $A_1\in\mathcal A_1$. Its
probability distribution (or, simply, distribution) $P^X$ with
respect to a probability measure $P$ on $\mathcal A$ is the image
measure of $P$ by $X$, i.e., the probability measure on $\mathcal
A_1$ defined by $P^X(A_1):=P(X^{-1}(A_1))$. Let us write $\times$
instead of $\otimes$ for the product of $\sigma$-fields or measures. The next definition is well known and can be found, for instance, in Heyer (1982).

\begin{mydef}\rm  (i) (Markov kernel) A Markov kernel
	$M_1:(\Omega,\mathcal A)\pt   (\Omega_1,\mathcal A_1)$ is a map $M_1:\Omega\times\mathcal A_1\rightarrow[0,1]$ such that: 
	a) $\forall \omega\in\Omega$, $M_1(\omega,\cdot)$ is a  probability
	measure on
	$\mathcal A_1$; b) $\forall A_1\in\mathcal A_1$, $M_1(\cdot,A_1)$ is an $\mathcal A$-measurable map.\\
  (ii) (Diagonal product of Markov kernels) The diagonal product
	$$M_1\times M_2:(\Omega,\mathcal A)\pt
	(\Omega_1\times\Omega_2,\mathcal A_1\times\mathcal A_2)$$ of two
	Markov kernels $M_1:(\Omega,\mathcal A)\pt   (\Omega_1,\mathcal A_1)$
	and $M_2:(\Omega,\mathcal A)\pt   (\Omega_2,\mathcal A_2)$ is defined
	as the only Markov kernel such that
	$$(M_1\times  M_2)(\omega,A_1\times A_2)=M_1(\omega,A_1)\cdot M_2(\omega,A_2),\quad A_i\in\mathcal A_i,
	i=1,2.
	$$
 (iii)
	(Image of a Markov kernel) The image (let us also call it {\it
		probability distribution}) of a Markov kernel $M_1:(\Omega,\mathcal
	A,P)\pt   (\Omega_1,\mathcal A_1)$ on a probability space is the
	probability measure  $P^{M_1}$ on $\mathcal A_1$ defined by
	$P^{M_1}(A_1):=\int_{\Omega}M_1(\omega,A_1)\,dP(\omega)$.
\end{mydef}

\begin{mydef}\rm  (Independence of Markov kernels, Nogales (2013a))
	Let $(\Omega,\mathcal A,P)$ be a probability space. Two Markov
	kernels $M_1:(\Omega,\mathcal A,P)\pt   (\Omega_1,\mathcal A_1)$ and
	$M_2:(\Omega,\mathcal A,P)\pt   (\Omega_2,\mathcal A_2)$ are said to
	be independent if $P^{M_1\times M_2}=P^{M_1}\times P^{M_2}$. We write $M_1\cind M_2$ {  (or $M_1\cind_P M_2$)}. 
\end{mydef}

Given two random variables $X_i:(\Omega,\mathcal
A,P)\rightarrow(\Omega_i,\mathcal A_i)$, $i=1,2$, the conditional
distribution of $X_2$ given $X_1$, when it exists, is a Markov
kernel $M_1:(\Omega_1,\mathcal A_1)\pt (\Omega_2,\mathcal A_2)$ such
that $P(X_1\in A_1,X_2\in
A_2)=\int_{A_1}M_1(\omega_1,A_2)dP^{X_1}(\omega_1)$, for all
$A_1\in\mathcal A_1$ and $A_2\in\mathcal A_2$. We write
$P^{X_2|X_1=\omega_1}(A_2):=M_1(\omega_1,A_2)$. Reciprocally, every
Markov kernel is also a conditional distribution, as it is noted in (2013b). This paper also introduces the next definition.

\begin{mydef}\rm  (Conditional distribution of a Markov kernel given another)
	Let  $M_1:(\Omega,\mathcal A,P)\pt   (\Omega_1,\mathcal A_1)$ and
	$M_2:(\Omega,\mathcal A,P)\pt   (\Omega_2,\mathcal A_2)$ be two Markov
	kernels over the same probability space. The conditional distribution
	$P^{M_1|M_2}$ of $M_1$ given $M_2$ is defined as a Markov kernel
	$L:(\Omega_2,\mathcal A_2)\pt   (\Omega_1,\mathcal A_1)$ such that,
	for every pair of events $A_1\in\mathcal A_1$ and $A_2\in\mathcal A_2$,
	\begin{gather*} 
			\int_\Omega
			M_1(\omega,A_1)M_2(\omega,A_2)dP(\omega)=\int_{A_2}L(\omega_2,A_1)dP^{M_2}(\omega_2).
	\end{gather*}
\end{mydef}

\begin{myrem} \rm An interesting problem in this context is  the existence of such conditional distributions, something that happens under well known regularity conditions on the involved measurable spaces, e.g. $(\Omega,\mathcal A)$, or the corresponding measurable space $(\Omega_i,\mathcal A_i)$, is a standard Borel space. This is the same for both random variables and Markov kernels (see Nogales (2013b)). In the rest of the paper we will assume this when necessary. $\Box$\end{myrem}

\section{Conditional Independence}\label{ci}

Let us recall the definition of conditional independence for random variables; we refer to Dawid (1979), for instance, where some basic properties are also given. 

\begin{mydef} \rm Let $X_i:(\Omega,\mathcal A,P)\rightarrow (\Omega_i,\mathcal A_i)$, $i=1,2,3,$ be arbitrary random variables $X_1$ and $X_2$ are said to be conditional independent given $X_3$, and we write $X_1\cind X_2|X_3$ {  (or $X_1\cind_P X_2|X_3$ to be more precise),} if
$$P^{(X_1,X_2)|X_3}=P^{X_1|X_3}\times P^{X_2|X_3},\quad P^{X_3}-c.s.
$$
\end{mydef} 

We are now ready for the main result of the paper.

\begin{mytheo}\label{teo2} \rm If $X_1$ and $X_2$ are conditional independent given  $X_3$, then $X_1$ and $X_2$ are independent if, and only if, the Markov kernels $P^{X_1|X_3}$ and $P^{X_2|X_3}$ are $P^{X_3}$-independent.\end{mytheo}

\begin{myrem}\rm (Some reformulations of the three propositions involved in the previous theorem) By definition, $X_1\cind X_2\mid X_3$ means that,
	for every $A_i\in\mathcal A_i$, $1\le i\le 3$,
	$$\int_{A_3}P^{(X_1,X_2)|X_3=\omega_3}(A_1\times A_2)dP^{X_3}(\omega_3)=
	\int_{A_3}P^{X_1|X_3=\omega_3}(A_1)\cdot P^{X_2|X_3=\omega_3}(A_2) dP^{X_3}(\omega_3).
	$$
This is equivalent to 
	$$E[(f_1\circ X_1)\cdot(f_2\circ X_2)|X_3]=E[f_1\circ X_1|X_3]\cdot E[f_2\circ X_2|X_3],\quad P^{X_3}-c.s.,
	$$
	for every bounded real random variables  $f_i:(\Omega_i,\mathcal A_i)\rightarrow\mathbb R$, $i=1,2$.
	
	In particular, $X_1\cind X_2$ is equivalent to
	$$E[(f_1\circ X_1)\cdot(f_2\circ X_2)]=E[f_1\circ X_1]\cdot E[f_2\circ X_2].
	$$
	
Finally, $P^{X_1|X_3}\cind_{P^{X_3}}P^{X_2|X_3}$ means that, for every $A_i\in\mathcal A_i$, $1\le i\le 2$, 
	\begin{gather*}\begin{split}
	\int_{\Omega_3}P^{X_1|X_3=\omega_3}(A_1)&\cdot P^{X_2|X_3=\omega_3}(A_2) dP^{X_3}(\omega_3)\\
	&=
	\int_{\Omega_3}P^{X_1|X_3=\omega_3}(A_1)dP^{X_3}(\omega_3)\cdot \int_{\Omega_3}P^{X_2|X_3=\omega_3}(A_2) dP^{X_3}(\omega_3),
	\end{split}\end{gather*}
	which is equivalent to 
	$$\int_{\Omega_3}E(f_1\circ X_1|X_3)\cdot E(f_2\circ X_2|X_3)dP^{X_3}=E(f_1\circ X_1)\cdot E(f_2\circ X_2).
	$$
	for every pair of functions $f_1,f_2$ as above. {  As $E(E(f_i\circ X_i|X_3))=E(f_i\circ X_i)$, this is equivalent to the uncorrelatedness of  the conditional expectations given $X_3$ of every pair of real bounded measurable functions of $X_1$ and $X_2$.} {  Seen in this way, the independence of these two conditional distributions has a degree of difficulty comparable to other conditions that appear in the literature cited in the bibliography; for instance, see  Proposition 2.4.g of van Putten et al. (1985), or others appearing in the results following it. }
%
\end{myrem}

\section{Counterexamples}\label{cexa}

Let $X_i:(\Omega,\mathcal A,P)\rightarrow (\Omega_i,\mathcal A_i)$, $i=1,2,3,$ be random variables. Consider the propositions:
	\begin{itemize}
		\item[(i)] $X_1\cind X_2\mid X_3$.
		\item[(ii)] $X_1\cind X_2$.
		\item[(iii)] $P^{X_1|X_3}\cind_{P^{X_3}}P^{X_2|X_3}$.
		\end{itemize}
		
		We have shown that $\mbox{(i)}+\mbox{(ii)}\Longrightarrow \mbox{(iii)}$ and $\mbox{(i)}+\mbox{(iii)}\Longrightarrow \mbox{(ii)}$, i.e., in presence of (i), the statements (ii) and (iii) are equivalent. In particular, (iii) is just we need to reach independence from conditional independence. 
		
		We can ask ourselves if every two of these propositions implies the third. In particular, we wonder if (i) and (ii) are equivalent when (iii) is satisfied. All the answers are negative, as the next counterexamples show. {  We also include two examples in which the theorem applies.}
		
		First, let us describe a common framework for them.
		
		Let $\Omega$ be a population with $n$ individuals and consider a partition $(A_{ijk})_{i,j,k=0,1}$ of $\Omega$. We write $n_{ijk}$ for the number of individuals of $A_{ijk}$. One or more of the indices $i,j,k$ can be replaced by a $+$ sign to denote the union of the corresponding sets of the partition: for instance, $A_{+01}=A_{001}\cup A_{101}$. In particular, $\Omega=A_{+++}$. Similar notations should be used for the numbers $n_{ijk}$ (e.g. $n_{+0+}=n_{000}+n_{001}+n_{100}+n_{101}$). Such a situation will be referred to as 
		$$C(n_{000},n_{001},n_{010},n_{011},n_{100},n_{101},n_{110},n_{111}).
		$$
		
		We introduce three dichotomic random variables $X_1,X_2,X_3$ as follows:
				\begin{gather*} 
				X_1(\omega)=i,\mbox{ if }\omega\in A_{i++}, i=0,1,\\
		X_2(\omega)=j,\mbox{ if }\omega\in A_{+j+}, j=0,1,\\
		X_3(\omega)=k,\mbox{ if }\omega\in A_{++k}, k=0,1.
		\end{gather*}
		
		\begin{myexa}\rm A scheme like this could be obtained when we want to study the relationship between two diagnostic procedures, represented by the dichotomous variables $X_1$ and $X_2$ ($X_i=1$ or $0$ when the $i^{th}$ diagnostic test is positive or negative, respectively), for a disease represented by the dichotomous variable $X_3$, which takes the values 1 or 0 depending on whether the disease is actually present or absent.  In this case, we have the following equivalence for some known related concepts:
			\begin{gather*}
			\mbox{prevalence of the disease}=\frac{n_{++1}}{n_{+++}},\\ 
			\mbox{specificity of }X_1=\frac{n_{+00}}{n_{++0}},\ \ \mbox{specificity of }X_2=\frac{n_{0+0}}{n_{++0}},\\ \mbox{sensitivity of }X_1=\frac{n_{+11}}{n_{++1}},\ \ \mbox{sensitivity of }X_2=\frac{n_{1+1}}{n_{++1}}.\quad \Box
			\end{gather*}	
			\end{myexa}
		
		The independence of $X_1$ and $X_2$ means that, for every $i,j=0,1$,
		$$n_{ij+}\cdot n_{+++}=n_{i++}\cdot n_{+j+}
		$$
		
The independence of $M_1:=P^{X_1|X_3}$ and $M_2:=P^{X_2|X_3}$ with respect to $P^{X_3}$ means that, for every $i,j=0,1$,
\begin{gather*}\sum_{k=0}^1P(X_1=i|X_3=k)\cdot P(X_2=j|X_3=k)\cdot P(X_3=k)=\\
\left(\sum_{k=0}^1P(X_1=i|X_3=k)\cdot P(X_3=k)\right)\cdot \left(\sum_{k=0}^1P(X_2=j|X_3=k)\cdot P(X_3=k)\right),
\end{gather*}
that is to say,
		\begin{gather*}
\frac{n_{0+0} n_{+00}}{n_{++0}}+\frac{n_{0+1} n_{+01}}{n_{++1}}=\frac{n_{0++} n_{+0+}}{n_{+++}}\\
\frac{n_{0+0} n_{+10}}{n_{++0}}+\frac{n_{0+1} n_{+11}}{n_{++1}}=\frac{n_{0++} n_{+1+}}{n_{+++}}\\
\frac{n_{1+0} n_{+00}}{n_{++0}}+\frac{n_{1+1} n_{+01}}{n_{++1}}=\frac{n_{1++} n_{+0+}}{n_{+++}}\\
\frac{n_{1+0} n_{+10}}{n_{++0}}+\frac{n_{1+1} n_{+11}}{n_{++1}}=\frac{n_{1++} n_{+1+}}{n_{+++}}
		\end{gather*}
		The conditional independence of $X_1$ and $X_2$ given $X_3$, i.e.  $P^{(X_1,X_2)|X_3}=P^{X_1|X_3}\times P^{X_2|X_3}$, means that, for every  $i,j,k=0,1$,
		$$P(X_1=i,X_2=j|X_3=k)=P(X_1=i|X_3=k)\cdot P(X_2=j|X_3=k),$$ or
		$$P(X_1=i,X_2=j,X_3=k)\cdot P(X_3=k)=P(X_1=i,X_3=k)\cdot P(X_2=j,X_3=k),
		$$
which is the same as
%
		{\begin{gather*}
		n_{000}\cdot n_{++0}=n_{0+0}\cdot n_{+00},\qquad n_{001}\cdot n_{++1}=n_{0+1}\cdot n_{+01}\\
		n_{010}\cdot n_{++0}=n_{0+0}\cdot n_{+10},\qquad n_{011}\cdot n_{++1}=n_{0+1}\cdot n_{+11}\\
		n_{100}\cdot n_{++0}=n_{1+0}\cdot n_{+00},\qquad n_{101}\cdot n_{++1}=n_{1+1}\cdot n_{+01}\\
		n_{110}\cdot n_{++0}=n_{1+0}\cdot n_{+10},\qquad n_{111}\cdot n_{++1}=n_{1+1}\cdot n_{+11}
	\end{gather*}}
	
The following counterexamples delimit Theorem \ref{teo2}. 
		
\begin{mycexa} \rm For $C(3000,200,1500,300,1500,200,3000,300)$ it is easy to see that $M_1=P^{X_1|X_3}$ and $M_2=P^{X_2|X_3}$ are  $P^{X_3}$-independent, but $X_1$ and $X_2$ are not $P$-independent. So, in absence of (i), (ii) is not implied by (iii). $\Box$
	\end{mycexa}

\begin{mycexa} \rm For $C(4200,400,2000,300,2000,200,1000,100)$,  $M_1=P^{X_1|X_3}$ and  $M_2=P^{X_2|X_3}$ are not  $P^{X_3}$-independent. Nevertheless $X_1$ and $X_2$ are independent. Obviously, $X_1$ and $X_2$ are not conditionally independent given $X_3$. So, in absence of (i), (iii) is not implied by (ii). $\Box$
\end{mycexa}

\begin{mycexa} \rm For $C(1000,1000,0,2000,0,2000,1000,1000)$,   $M_1=P^{X_1|X_3}$ and $M_2=P^{X_2|X_3}$  are $P^{X_3}$-independent, and $X_1$ and $X_2$ are independent, but $X_1$ and $X_2$ are not conditionally independent given $X_3$. So (i) is not implied by (ii)+(iii). $\Box$
\end{mycexa}

{  In the next two examples, the condition (i) holds. Hence the propositions (ii) and (iii) hold or not simultaneously. See also the remark below to see how Theorem 1 is an improvement of two previous known results on the relationship between unconditional and conditional independence.}

{
\begin{myexa}\rm  For $C(1200,3000,1200,3000,2000,3200,2000,3200)$  the three propositions (i), (ii) and (iii) are satisfied.  $\Box$
\end{myexa}}

{
\begin{myexa}\rm  For $C(1200,3000,1200,3200,2000,3000,2000,3200)$, (i) holds, but not  (ii) or (iii).  $\Box$
\end{myexa}}

\begin{myrem}\rm   Keeping the previous notations, it is known that (i) + $X_1\cind X_3$ implies (ii); see, for instance, Florens et al. (2000, Theorem 2.2.10) or Lemma 4.3 of Dawid (1979) when the conditioning on $Z$ is absent. Theorem \ref{teo2} is an improvement of this result as $X_1\cind X_3$ implies, and it is not implied {  by,} (iii), as we prove in what follows.	It is easy to see that the independence of $X_1$ and $X_3$ implies (iii). Indeed, given bounded real random variables $f_i$, $i=1,2$, the independence of $X_1$ and $X_3$ yields $E(f_1\circ X_1|X_3)=E(f_1\circ X_1)$ and hence 
$$\int_{\Omega_3}E(f_1\circ X_1|X_3)E(f_2\circ X_2|X_3)dP^{X_3} 
=E(f_1\circ X_1)E(f_2\circ X_2),$$ which is equivalent to (iii). Let us show that the reciproque is not true: it is proved in Nogales (2013b) that, for a trivariate normal random variable $(X_1,X_2,X_3)$ with null mean and covariance matrix $(\sigma_{ij})$, the $P^{X_3}$-conditional distribution $L(x_1,\cdot)$ of $P^{X_2|X_3}$ given that $P^{X_1|X_3}$ has taken the value $x_1$ follows a normal distribution with mean $\sigma_1^{-1}\sigma_2\rho_{23}\rho_{13}x_1$ and variance $\sigma_2^2(1-\rho_{23}^2\rho_{13}^2)$, {  where $\rho_{ij}$ denotes the correlation coefficient of $X_i$ and $X_j$}. So, the Markov kernels $P^{X_2|X_3}$ and $P^{X_1|X_3}$ are $P^{X_3}$-independent if, and only if, $L(x_1,\cdot)$ coincides with $(P^{X_3})^{P^{X_2|X_3}}$ (which coincides with  $P^{X_2}$), and this happens if $\rho_{23}=0$ or $\rho_{13}=0$. So, for $\rho_{23}=0$ and $\rho_{13}\ne 0$, we have that  $P^{X_2|X_3}$ and $P^{X_1|X_3}$ are $P^{X_3}$-independent, but $X_1$ and $X_3$ are not independent. $\Box$
	\end{myrem}

 \begin{myrem} \rm  Phillips (1988) shows the next result: ``For $i=1,2$, consider random variables 	$Y_i:(\Omega,\mathcal A,P)\rightarrow (\Omega_i,\mathcal A_i)$,  $Z_i:(\Omega_1\times\Omega_2,\mathcal A_1\times\mathcal A_2)\rightarrow(\Omega'_i,\mathcal A'_i)$, $f_i:(\Omega'_i,\mathcal A'_i)\rightarrow(\Omega''_i,\mathcal A''_i).$ If $f_1\circ Z_1\circ(Y_1,Y_2)\cind f_2\circ Z_2\circ(Y_1,Y_2)|Y_1$, then 
	$$f_1\circ Z_1\circ(Y_1,Y_2)\cind f_2\circ Z_2\circ(Y_1,Y_2)$$ is equivalent to
	$$E(E(I_{F_1}|Y_1)I_{F_2})=E(I_{F_1})\cdot E(I_{F_2})\quad \text{and}\quad E(E(I_{F_2}|Y_1)I_{F_1})=E(I_{F_1})\cdot E(I_{F_2}),
	$$
	whatever be the events $F_i\in (f_i\circ Z_i\circ(Y_1,Y_2))^{-1}(\mathcal A''_i)$, $i=1,2.$''
	
	This is a particular case of Theorem \ref{teo2} with no more to take $X_1=f_1\circ Z_1\circ(Y_1,Y_2)$, $X_2=f_2\circ Z_2\circ(Y_1,Y_2)$ y $X_3=Y_1$. Indeed, according to Theorem \ref{teo2}, if $X_1\cind X_2|X_3$, then $X_1\cind X_2$ is equivalent to $P^{X_1|X_3}\cind_{P^{X_3}} P^{X_2|X_3}$, which in turns means that, for every bounded real random variable $g_i$ on $(\Omega''_i,\mathcal A''_i)$, 
	$$\int_{\Omega_3}E(g_1\circ X_1|X_3)\cdot E(g_2\circ X_2|X_3)dP^{X_3}=E(g_1\circ X_1)\cdot E(g_2\circ X_2).
	$$
	If $F_i=X_i^{-1}(A''_i)$, making $g_i:=I_{A''_i}$, $i=1,2$, it follows that 
	$$E[E(I_{F_1}|X_1)E(I_{F_2}|X_1)]=E(I_{F_1})\cdot E(I_{F_2}),
	$$
	and, on the other hand,
	$$E[E(I_{F_1}|X_1)E(I_{F_2}|X_1)]=E[E(I_{F_1}E(I_{F_2}|X_1)|X_1)]=E[I_{F_1}E(I_{F_2}|X_1)]
	$$
	and
	$$E[E(I_{F_1}|X_1)E(I_{F_2}|X_1)]=E[E(I_{F_2}E(I_{F_1}|X_1)|X_1)]=E[I_{F_2}E(I_{F_1}|X_1)].
	$$
Hence
	$$E[I_{F_1}E(I_{F_2}|X_1)]=E[I_{F_2}E(I_{F_1}|X_1)]=E(I_{F_1})\cdot E(I_{F_2}).$$
	It is readily shown that, from these two equalities, we obtain   
	$$\int_{\Omega_3}E(g_1\circ X_1|X_3)\cdot E(g_2\circ X_2|X_3)dP^{X_3}=E(g_1\circ X_1)\cdot E(g_2\circ X_2),
	$$
for every bounded real random variables $g_1,g_2$ on $(\Omega''_i,\mathcal A''_i)$. $\Box$
\end{myrem}

\section{Extension to Markov kernels}\label{extmk}

In this section we extend to Markov kernels the concept of conditional independence. Theorem \ref{teo2} is also extended to this framework. 

\begin{mydef}\rm (Conditional independence of Markov kernels) Given three Markov kernels $M_i:(\Omega,\mathcal A,P)\pt (\Omega_i,\mathcal A_i)$, $1\le i\le 3$, we shall say that $M_1$ and $M_2$ are conditionally independent given $M_3$, and we write $M_1\cind_P M_2|M_3$ {  (or $M_1\cind M_2|M_3$ if there is not ambiguity),}  when
$$P^{M_1\times M_2|M_3}=P^{M_1|M_3}\times P^{M_2|M_3},\quad P^{M_3}-c.s.
$$
\end{mydef}

\begin{myrem}\rm (A representation in terms of random variables) Keeping the suppositions of the previous definition, let us write $q_i$ for the natural $i^{\mbox{th}}$ projection on $\Omega_1\times\Omega_2\times\Omega_3$, $1\le i\le 3.$ It is readily shown that 
\begin{gather*}
\big(P^{M_1\times M_2\times M_3}\big)^{q_i}=P^{M_i},\quad 1\le i\le 3,\vspace{1.5ex}\\
P^{M_1\times M_2|M_3}=\big(P^{M_1\times M_2\times M_3}\big)^{(q_1,q_2)|q_3},\vspace{1.5ex}\\
P^{M_i|M_3}=\big(P^{M_1\times M_2\times M_3}\big)^{q_i|q_3},\quad i=1,2.
\end{gather*}
So,
$$M_1\cind_P M_2|M_3\quad\Longleftrightarrow\quad q_1\cind_{P^{M_1\times M_2\times M_3}} q_2|q_3.
$$

Moreover, when $\Omega_2=\mathbb R^k$ and $M_2$ is integrable, from
$$E(M_i|M_1)(\omega_1)=\int_{\mathbb R^k} xdP^{M_i|M_1=\omega_1}(x)
$$
we obtain that
$$E_P(M_i|M_1)=E_{P^{M_1\times M_2\times M_3}}(q_i|q_1),\quad i=1,2. \   \ \Box
$$
	\end{myrem}

\begin{myrem}\rm (Characterization in terms of densities) Suppose that, for $i=1,2,3$,  $\mu_i$ is a $\sigma$-finite measure on
	$\mathcal A_i$ such that
	$dM_i(\omega,\cdot)=\phi_i(\omega,\cdot)d\mu_i$, where $\phi_i$ is a
	nonnegative real $\mathcal A\times\mathcal A_i$-measurable function
	on $\Omega\times\Omega_i$. Usually, the dominating measure $\mu_i$ is the counting measure in the discrete (respectively, the Lebesgue measure in the continuous) case, both in the univariate and multivariate framework. It is shown in Nogales (2013b) that the map $\omega_3\mapsto\int_\Omega
	\phi_3(\omega,\omega_3)dP(\omega)$ is a $\mu_3$-density of
	$P^{M_3}$ and, besides, for $i=1,2$, the conditional distribution
	$L_i:=P^{M_i|M_3}$ exists and, for $P^{M_3}$-almost every $\omega_3$, the map
	$$\omega_i\mapsto
	\frac{\int_\Omega\phi_i(\omega,\omega_i)\phi_3(\omega,\omega_3)dP(\omega)}{\int_\Omega\phi_3(\omega,\omega_3)dP(\omega)}
	$$
	is a $\mu_i$-density of $L_i(\omega_3,\cdot)$.
	
	A similar reasoning shows that the map
	$$(\omega_1,\omega_2)\mapsto \int_\Omega \phi_1(\omega,\omega_1)\phi_2(\omega,\omega_2)dP(\omega)
	$$
	is a $\mu_1\times\mu_2$-density of $P^{M_1\times M_2}$, and the conditional distribution
	$L:=P^{M_1\times M_2|M_3}$ exists and, for $P^{M_3}$-almost every $\omega_3$, the map
	$$(\omega_1,\omega_2)\mapsto
	\frac{\int_\Omega\phi_1(\omega,\omega_1)\phi_2(\omega,\omega_2)\phi_3(\omega,\omega_3)dP(\omega)}{\int_\Omega\phi_3(\omega,\omega_3)dP(\omega)}
	$$
	is a $\mu_1\times\mu_2$-density of $L(\omega_3,\cdot)$.\par
Hence, the conditional independence of $M_1$ and $M_2$ given $M_3$ means that, for $P^{M_3}$-almost every $\omega_3$ and $\mu_1\times\mu_2$-almost every $(\omega_1,\omega_2)$,
	$$\frac{\int_\Omega\phi_1(\omega,\omega_1)\phi_2(\omega,\omega_2)\phi_3(\omega,\omega_3)dP(\omega)}{\int_\Omega\phi_3(\omega,\omega_3)dP(\omega)}=\frac{\int_\Omega\phi_1(\omega,\omega_1)\phi_3(\omega,\omega_3)dP(\omega)}{\int_\Omega\phi_3(\omega,\omega_3)dP(\omega)}\cdot \frac{\int_\Omega\phi_2(\omega,\omega_2)\phi_3(\omega,\omega_3)dP(\omega)}{\int_\Omega\phi_3(\omega,\omega_3)dP(\omega)}.\ \Box
	$$
	\end{myrem}

The next theorem extends Theorem \ref{teo2} to Markov kernels.

\begin{mytheo}\label{teo3}\rm Let $M_i:(\Omega,\mathcal A,P)\pt (\Omega_i,\mathcal A_i)$, $i=1,2,3,$ be Markov kernels. Consider the propositions:
	\begin{itemize}
		\item[(i)] $M_1\cind M_2\mid M_3$.
		\item[(ii)] $M_1\cind M_2$.
		\item[(iii)] $P^{M_1|M_3}\cind_{P^{M_3}}P^{M_2|M_3}$.
	\end{itemize}
Then, under (i), the propositions (ii) and (iii) are equivalent.
\end{mytheo}

\section{Another extension of the main result}\label{last}

 A more general result than Theorem \ref{teo2} in terms of random variables is presented in this section, where the introduced definition of conditional independence between Markov kernels is used to obtain a minimal condition which added to conditional independence of $X_1$ and $X_2$ given $X_3$ implies the conditional independence of $X_1$ and $X_2$ given $X_4$ when $X_4$ is function of $X_3$. In fact,  Theorem \ref{teo2} appears as the particular case  in which  $X_4$ is a constant function.

\begin{mytheo}\label{teo4}\rm Let $X_i:(\Omega,\mathcal A,P)\pt (\Omega_i,\mathcal A_i)$, $i=1,2,3,4,$ random variables. Suppose that $X_4=f\circ X_3$, where  $f:(\Omega_3,\mathcal A_3)\rightarrow(\Omega_4,\mathcal A_4)$. Consider the propositions:
	\begin{itemize}
		\item[(i)] $X_1\cind X_2\mid X_3$.
		\item[(ii)] $X_1\cind X_2\mid X_4$.
		\item[(iii)] $P^{X_1|X_3}\cind_{P^{X_3}}P^{X_2|X_3}\mid P^{X_4|X_3}$.
	\end{itemize}
	Then, if (i) holds, the propositions (ii) and (iii) are  equivalent. 
\end{mytheo}

\begin{myrem}\rm To obtain a characterization of the statement (iii), note first that, for $i=1,2$,
	$$(P^{X_3})^{P^{X_i|X_3}|P^{X_4|X_3}}(\cdot,A_i)=E[E(I_{A_i}\circ X_i|\sigma(X_3))|X_4],
	$$
	where $\sigma(X_3)$ denotes the $\sigma$-field $X_3^{-1}(\mathcal A_3)$ induced by $X_3$. 
	Indeed, we have that, by definition,  $(P^{X_3})^{P^{X_i|X_3}|P^{X_4|X_3}}$ ($:=Q^{M_i|M_4}$) is a Markov kernel  $M_{i4}:(\Omega_4,\mathcal A_4)\pt (\Omega_i,\mathcal A_i)$ such that
	$$\int_{\Omega_3}P^{X_i|X_3}(\omega_3,A_i)\cdot P^{X_4|X_3}(\omega_3,A_4)dP^{X_3}(\omega_3)=
	\int_{A_4}M_{i4}(\omega_4,A_i)dQ^{M_4}(\omega_4),
	$$ 
	for every $A_i\in\mathcal A_i$ and $A_4\in\mathcal A_4.$
	But, 
	\begin{gather*}\begin{split}
	\int_{\Omega_3}P^{X_i|X_3}(\omega_3,A_i)&\cdot P^{X_4|X_3}(\omega_3,A_4)dP^{X_3}(\omega_3)\\
	&=\int_{\Omega_3} E(I_{A_i}\circ X_i|X_3)\cdot E(I_{A_4}\circ X_4|X_3)dP^{X_3}\\
	&=\int_{\Omega} E(I_{A_i}\circ X_i|\sigma(X_3))\cdot E(I_{A_4}\circ X_4|\sigma(X_3))dP\\
	&=\int_{\Omega}  E\big[(I_{A_4}\circ X_4)\cdot E(I_{A_i}\circ X_i|\sigma(X_3))\big|\sigma(X_3)\big]dP\\
	&=\int_{X_4^{-1}(A_4)}  E(I_{A_i}\circ X_i|\sigma(X_3))dP.
	\end{split}\end{gather*}

	Since $Q^{M_4}=P^{X_4}$, it readily follows that 
	$$(P^{X_3})^{P^{X_i|X_3}|P^{X_4|X_3}}(\cdot,A_i)=E[E(I_{A_i}\circ X_i|\sigma(X_3))|X_4].
	$$
Analogously, by definition, 
$(P^{X_3})^{P^{X_1|X_3}\times P^{X_2|X_3}|P^{X_4|X_3}}$ ($:=Q^{M_1\times M_2|M_4}$) is a Markov kernel  $M_{(12)4}:(\Omega_4,\mathcal A_4)\pt (\Omega_1\times \Omega_2,\mathcal A_1\times \mathcal A_1)$ such that
$$\int_{\Omega_3}P^{X_1|X_3}(\omega_3,A_1)\cdot P^{X_2|X_3}(\omega_3,A_2)\cdot P^{X_4|X_3}(\omega_3,A_4)dP^{X_3}(\omega_3)=
\int_{A_4}M_{(12)4}(\omega_4,A_1\times A_2)dQ^{M_4}(\omega_4),
$$ 
for every $A_i\in\mathcal A_i$ , $i=1,2,4$. But
%
	\begin{gather*}\begin{split}
		\int_{\Omega_3} P^{X_1|X_3}(\omega_3,A_1)&\cdot P^{X_2|X_3}(\omega_3,A_2)\cdot P^{X_4|X_3}(\omega_3,A_4)dP^{X_3}(\omega_3)\\&=\int_{\Omega_3} E(I_{A_1}\circ X_1|X_3)\cdot E(I_{A_2}\circ X_2|X_3)\cdot E(I_{A_4}\circ X_4|X_3)dP^{X_3}\\
		&=\int_{\Omega} E(I_{A_1}\circ X_1|\sigma(X_3))\cdot E(I_{A_2}\circ X_2|\sigma(X_3))\cdot E(I_{A_4}\circ X_4|\sigma(X_3))dP\\
		&=\int_{\Omega}  E\big[(I_{A_4}\circ X_4)\cdot E(I_{A_1}\circ X_1|\sigma(X_3))\cdot E(I_{A_2}\circ X_2|\sigma(X_3))\big|\sigma(X_3)\big]dP\\
		&=\int_{X_4^{-1}(A_4)}  E(I_{A_1}\circ X_1|\sigma(X_3))\cdot E(I_{A_2}\circ X_2|\sigma(X_3))dP.
\end{split}\end{gather*}

	So, the statement (iii) $P^{X_1|X_3}\cind_{P^{X_3}}P^{X_2|X_3}\mid P^{X_4|X_3}$ can be expressed in the form
	$$E[E(f_1\circ X_1|\sigma(X_3))\cdot E(f_2\circ X_2|\sigma(X_3))|X_4]=E[E(f_1\circ X_1|\sigma(X_3))|X_4]\cdot E[E(f_2\circ X_2|\sigma(X_3))|X_4], 
	$$
	for every bounded real random variable  $f_i$ on $(\Omega_i,\mathcal A_i)$, $i=1,2$. 
	$\Box$
\end{myrem}

{  In the previous result the $\sigma$-field $\sigma(X_4)$ is contained in  $\sigma(X_3)$. The reader is referred to {van Putten} et al. (1985) where invariance properties of conditional independence under enlargement or reduction is systematically investigated.}

\section{Proofs}

{\sc Proof of Theorem \ref{teo2}.} \rm  Let us write $Q=P^{X_3}$ and $M_i=P^{X_i|X_3}$, $i=1,2$.	In the following, we suppose  $X_1\cind X_2|X_3$.

1) We show first that if $X_1\cind X_2$,  then 
$$Q^{M_1\times M_2}=Q^{M_1}\times Q^{M_2}.
$$
Note that 
$$Q^{M_i}(A_i)=\int_{\Omega_3}M_i(\omega_3,A_i)dQ(\omega_3)=
\int_{\Omega_3}P^{X_i|X_3=\omega_3}(A_i)dP^{X_3}(\omega_3)=
P^{X_i}(A_i),$$
i.e., 
$$\left(P^{X_3}\right)^{P^{X_i|X_3}}=P^{X_i}.
$$
By definition,
$$(M_1\times M_2)(\omega_3,A_1\times A_2)=M_1(\omega_3,A_1)\cdot 
M_2(\omega_3,A_2).
$$
Hence, by conditional independence, 
\begin{gather*}\begin{split}
Q^{M_1\times M_2}(A_1\times A_2)&=
\int_{\Omega_3}P^{X_1|X_3=\omega_3}(A_1)\cdot 
P^{X_2|X_3=\omega_3}(A_2)dP^{X_3}(\omega_3)\\
&=\int_{\Omega_3}P^{(X_1,X_2)|X_3=\omega_3}(A_1\times A_2)dP^{X_3}(\omega_3)\\
&=P^{(X_1,X_2)}(A_1\times A_2),
\end{split}
\end{gather*}
which coincides with
$$P^{X_1}(A_1)\cdot P^{X_2}(A_2)=\left(P^{X_3}\right)^{P^{X_1|X_3}}(A_1)\cdot \left(P^{X_3}\right)^{P^{X_2|X_3}}(A_2)=Q^{M_1}(A_1)\cdot Q^{M_2}(A_2)
$$
since  $X_1$ and $X_2$ are independent. 

2) Now suppose that  the Markov kernels $P^{X_1|X_3}$ and $P^{X_2|X_3}$ are $P^{X_3}$-independent (in addition that $X_1\cind X_2|X_3$). Then, given $A_i\in\mathcal A_i$, $i=1,2$, we have that
\begin{gather*}\begin{split}
		P^{X_1}(A_1)\cdot P^{X_2}(A_2)&=Q^{M_1}(A_1)\cdot Q^{M_2}(A_2)\\
		&=Q^{M_1\times M_2}(A_1\times A_2)\\
		&=\int_{\Omega_3}P^{X_1|X_3=\omega_3}(A_1)\cdot 
		P^{X_2|X_3=\omega_3}(A_2)dP^{X_3}(\omega_3)\\
		&=\int_{\Omega_3}P^{(X_1,X_2)|X_3=\omega_3}(A_1\times A_2)dP^{X_3}(\omega_3)\\
		&=P^{(X_1,X_2)}(A_1\times A_2),
	\end{split}
\end{gather*}
which shows that $X_1\cind X_2$.
$\Box$\vspace{2ex}

{\sc Proof of Theorem \ref{teo3}.} \rm Let $q_i:\Omega_1\times\Omega_2\times\Omega_3\rightarrow \Omega_i$ the natural $i^{\text{\footnotesize th}}$ projection, $1\le i \le 3$. Writing $Q=
	P^{M_1\times M_2\times M_3}$, we have that
\begin{gather*}
Q^{q_i}=P^{M_i},\quad 1\le i\le 3,\\
Q^{(q_1,q_2)}=P^{M_1\times M_2},\\
P^{M_1\times M_2|M_3}=Q^{(q_1,q_2)|q_3},\\
P^{M_i|M_3}=Q^{q_i|q_3},\quad i=1,2.
\end{gather*}
It follows that
\begin{gather*}
	M_1\cind_P M_2|M_3\quad\Longleftrightarrow\quad q_1\cind_Q q_2|q_3,\\
	M_1\cind_P M_2\quad\Longleftrightarrow\quad q_1\cind_Q q_2,\\
	P^{M_1|M_3}\cind_{P^{M_3}}P^{M_2|M_3}\quad\Longleftrightarrow\quad Q^{q_1|q_3}\cind_{Q^{q_3}}Q^{q_2|q_3},
	\end{gather*}
and the result becomes a consequence of this and Theorem \ref{teo2}. $\Box$
\vspace{2ex}

{\sc Proof of Theorem \ref{teo4}.} \rm  	
	Consider the Markov kernels $M_i=P^{X_i|X_3}:(\Omega_3,\mathcal A_3)\pt (\Omega_i,\mathcal A_i)$, $i=1,2,4.$ Write $Q=P^{X_3}$. 
	Note that proposition (iii), i.e. $M_1\cind_Q M_2\mid M_4$, means that
	$$Q^{M_1\times M_2|M_4}=Q^{M_1|M_4}\times Q^{ M_2|M_4}.
	$$
We assume that (i) holds, that is, $P^{(X_1,X_2)|X_3}=P^{X_1|X_3}\times P^{X_2|X_3}$. Under such assumption, it will be enough to prove that
$$Q^{M_1\times M_2|M_4}=P^{(X_1,X_2)|X_4},\quad\text{and}\quad Q^{M_i|M_4}=P^{X_i|X_4},\ \ i=1,2.
$$
Let us show the first equality, the second being similar. 

By definition, 
$Q^{M_1\times M_2|M_4}$ is a Markov kernel $M:(\Omega_4,\mathcal A_4)\pt (\Omega_1\times\Omega_2,\mathcal A_1\times\mathcal A_2)$ such that, for every $C\in\mathcal A_1\times\mathcal A_2$ and $A_4\in\mathcal A_4,$
	$$\int_{\Omega_3}(M_1\times M_2)(\omega_3,C)\cdot M_4(\omega_3,A_4)dQ(\omega_3)=
\int_{A_4}M(\omega_4,C)dQ^{M_4}(\omega_4).
$$
Note that, as it can be easily verified, $Q^{M_4}=P^{X_4}$. Note also that, being $X_4=f\circ X_3$, $$M_4(\omega_3,A_4)=P^{X_4|X_3}(\omega_3,A_4)=I_{f^{-1}(A_4)}(\omega_3),\quad Q-a.s.$$
So, 
$$\int_{\Omega_3}(M_1\times M_2)(\omega_3,C)\cdot M_4(\omega_3,A_4)dQ(\omega_3)=\int_{f^{-1}(A_4)}(M_1\times M_2)(\omega_3,C)dP^{X_3}(\omega_3).
$$

	It follows that
	$$\int_{f^{-1}(A_4)}(M_1\times M_2)(\omega_3,C)dP^{X_3}(\omega_3)=\int_{A_4}M(\omega_4,C)dP^{X_4}(\omega_4).
	$$
	Moreover, using (i),
	\begin{gather*}\begin{split}\int_{f^{-1}(A_4)}(M_1\times M_2)(\omega_3,C)dP^{X_3}(\omega_3)&=\int_{f^{-1}(A_4)}\big(P^{X_2|X_3}\times P^{X_2|X_3}\big)(\omega_3,C)dP^{X_3}(\omega_3)\\
	&=\int_{f^{-1}(A_4)}P^{(X_1,X_2)|X_3}(\omega_3,C)dP^{X_3}(\omega_3)\\
	&=P^{(X_1,X_2,X_3)}(C\times f^{-1}(A_4))=P^{(X_1, X_2,X_4)}(C\times A_4)\\
	&=\int_{A_4}P^{(X_1,X_2)|X_4}(\omega_4,C)dP^{X_4}(\omega_4),
	\end{split}\end{gather*}
	which shows that, $Q^{M_1\times M_2|M_4}=P^{(X_1,X_2)|X_4}.$\par
An analogous reasoning ((i) is not needed in this case) shows that $Q^{M_i|M_4}=P^{X_i|X_4},\ \ i=1,2,$ and this finishes the proof. $\Box$
\vspace{2ex}

\section{Acknowledgments}
This work was supported by the  \emph{Junta de Extremadura} (Autonomous Government of Extremadura, Spain) under the project GR15013.

\section* {References:}

\begin{itemize}

\item Dawid, A.P. (1979) Conditional Independence in Statistical Theory, Journal of the Royal Statistical Society B 41, 1-31.

\item Dawid, A.P. (1980) Conditional Independence for Statistical Operations, Annals of Statistics 8, 598-617.

\item Florens, J.P., Mouchart, M., and Rolin, J.M. (1990) Elements of Bayesian Statistics, Marcel Dekker, New York.

	\item Heyer, H. (1982)  Theory of Statistical Experiments, Springer, Berlin.

	\item Nogales, A.G.  (2013a) On Independence of Markov Kernels and a Generalization of Two Theorems of Basu, Journal of Statistical Planning and Inference 143, 603-610.

	\item Nogales, A.G. (2013b) Existence of Regular Conditional Probabilities for Markov Kernels, Statistics and Probability Letters 83, 891-897.
	
	\item Phillips, P.C.B. (1988) Conditional and Unconditional Statistical Independence, Journal of Econometrics 38, 341-348.
	
	\item {  van Putten, C.; van Schuppen, J.H. {(1985)} Invariance Properties of the Conditional Independence Relation, Ann. Probab. 13, no. 3, 934--945.}

\end{itemize}

\end{document}